\documentclass[12pt]{article}

\begin{document}

\begin{titlepage}

\vspace*{-2cm}

\vspace{.5cm}

\begin{centering}

\huge{Some sufficient conditions of a given series with rational
terms converging to an irrational number or a transcdental number}

\vspace{.5cm}

\large  {Yun Gao,Jining Gao }\\

\vspace{.5cm}

Shanghai Putuo college, Shanghai Jiaotong University

\vspace{.5cm}

\begin{abstract}

In this paper, we propose various sufficient conditions to determine
if a given real number is an irrational number or a transcendental
number and also apply these conditions to some interesting examples
,particularly,one of them comes from complex analytic dynamics.
\end{abstract}

\end{centering}

\end{titlepage}

\pagebreak

\def\lh{\hbox to 15pt{\vbox{\vskip 6pt\hrule width 6.5pt height 1pt}
  \kern -4.0pt\vrule height 8pt width 1pt\hfil}}
\def\blob{\mbox{$\;\Box$}}
\def\qed{\hbox{${\vcenter{\vbox{\hrule height 0.4pt\hbox{\vrule width
0.4pt height 6pt \kern5pt\vrule width 0.4pt}\hrule height
0.4pt}}}$}}

\newtheorem{theorem}{Theorem}
\newtheorem{lemma}[theorem]{Lemma}
\newtheorem{definition}[theorem]{Definition}
\newtheorem{corollary}[theorem]{Corollary}
\newtheorem{proposition}[theorem]{Proposition}
\newcommand{\proof}{\bf Proof.\rm}

\section{Introduction}
In series theory, there is well known Cauchy convergence test which
is used to determine convergence of a given series,but Cauchy
convergence test usually is not practical in most applications,so
there come out various convergence test such as D'Alembert
convergence test ,integral convergence test and so on. In
Diophantine approximation theory, we are in totally different
situation that we already have necessary and sufficient condition
to determine if a given real number is an irrational number or a
transcendental number such as well known Roth theorem but seems to
be lack of practical test just as various convenient test in series
theory.
\newline
The purpose of this paper is to propose some  sufficient conditions
for convenient use in determining if a given real number is an
irrational number or a transcendental number and also give out
various interesting examples  to illustrate how to apply these
conditions,particularly, we will explain an example coming from
complex analytic dynamics in detail. At the end of this paper, we
propose a conjecture about rational approximation of any irrational
number.

\begin{theorem}
Assume that series $\sum_{n=1}^{\infty}c_{n},
c_{n}=\frac{a_{n}}{b_{n}}\neq 0 ,(n=1,2,\cdots)$ are rational
numbers and satisfy following conditions:
\newline
(1) $$b_n |b_{n+1},(n=1,2,\cdots)$$
\newline
(2)$$lim_{n\rightarrow\infty}a_n\frac{c_{n+1}}{c_{n}}=0$$
\newline
(3)for any natural number $ N$, $\sum_{n=N}^{\infty}c_{n}\neq 0$
\newline
then the series$\sum_{n=N}^{\infty}c_{n}$ converges to  an
irrational number.

\end{theorem}
\begin{proof} First of all, since condition (2) implies
$lim_{n\rightarrow\infty}\frac{c_{n+1}}{c_{n}}=0$, the series
$\sum_{n=1}^{\infty}c_{n}$ is convergent and set the convergent
result to be $\theta$. We will use indirect method to show $\theta$
is an irrational number. Suppose that $\theta=\frac{s}{r}$ is a
rational number. By (1), when $k\leq n$,$b_k |b_{n}$, we have
\begin{eqnarray}
A_n =
rb_{n}(\frac{s}{r}-\frac{a_1}{b_1}-\cdots-\frac{a_n}{a_n})\nonumber
\\ =rb_{n}(\frac{a_{n+1}}{b_{n+1}}+\frac{a_{n+2}}{b_{n+2}}+\cdots
)\nonumber \\=rb_{n}(c_{n+1}+c_{n+2}+\cdots)\neq 0\nonumber
\end{eqnarray}

and $A_n$ is an integer number, we notice
that$$A_n=rb_{n}c_{n}\frac{c_{n+1}}{c_{n}}(1+\frac{c_{n+2}}{c_{n+1}}+\frac{c_{n+3}}{c_{n+1}}+\cdots
)$$ and $lim_{n\rightarrow\infty}\frac{c_{n+1}}{c_{n}}=0$, thus
there exists $N_1$, when $n\geq N_1$ we have
$$\frac{|c_{n+1}|}{|c_{n}|}<\frac{1}{2},\frac{|c_{n+2}|}{|c_{n+1}|}<\frac{1}{2},\cdots$$,
so
$$\frac{|c_{n+3}|}{|c_{n+1}|}=\frac{|c_{n+3}|}{|c_{n+2}|}\frac{|c_{n+2}|}{|c_{n+1}|}<(\frac{1}{2})^2$$
$$\frac{|c_{n+4}|}{|c_{n+1}|}=\frac{|c_{n+4}|}{|c_{n+3}|}\frac{|c_{n+3}|}{|c_{n+2}|}\frac{|c_{n+2}|}{|c_{n+1}|}<(\frac{1}{2})^3$$

Therefore,$$|A_n|\leq |r||a_{n}|
\frac{c_{|n+1}|}{c_{|n|}}(1+\frac{c_{|n+2|}}{c_{|n+1|}}+\frac{c_{|n+3|}}{c_{|n+1|}}+\cdots
)$$
$$<|r||a_{n}|
\frac{c_{|n+1}|}{c_{|n|}}(1+\frac{1}{2}+\frac{1}{2^2}+\cdots )$$
$$=2|r||a_{n}|
\frac{c_{|n+1}|}{c_{|n|}}$$ By (2),for $\frac{1}{2|r|}>0$,there
exists $N_2$,when $n\geq N_2$,$|a_{n}|
\frac{c_{|n+1}|}{c_{|n|}}<\frac{1}{2|r|}$. Set $N=max(N_1,N_2)$,when
$n\geq N$,we get $|A_n|<2|r||a_{n}|
\frac{c_{|n+1}|}{c_{|n|}}<2r\frac{1}{2r}=1$ which contradicts the
fact that $A_n$ is an integer and $A_n\neq 0$. That follows the
theorem.
\end{proof}
When the series just contains positive terms, condition (3) is
naturally satisfied, we have

\begin{theorem}
Assume that series $\sum_{n=1}^{\infty}c_{n},
c_{n}=\frac{a_{n}}{b_{n}}> 0 ,(n=1,2,\cdots)$ are rational numbers
and satisfy following conditions:
\newline
 (1) $$b_n
|b_{n+1},(n=1,2,\cdots)$$
\newline
(2) $$lim_{n\rightarrow\infty}a_n\frac{c_{n+1}}{c_{n}}=0$$  then the
series$\sum_{n=1}^{\infty}c_{n}$ converges to an irrational number.
\end{theorem}

{\bf Remark.} In the above theorem, the condition
$lim_{n\rightarrow\infty}a_n\frac{c_{n+1}}{c_{n}}=0$ is not
sufficient because $\sum_{n=0}^{\infty}\frac{n+1}{n!}=2e$ is an
irrational number and
$a_n\frac{c_{n+1}}{c_{n}}=\frac{n+2}{n+1}\rightarrow 1(n\rightarrow
\infty)$

\bigskip

{\bf Example 1.}
$$e=\sum_{n=0}^{\infty}\frac{1}{n!}=1+\frac{1}{1!}+\frac{1}{2!}+\frac{1}{3!}+\cdots
+\frac{1}{n!}+\cdots$$ is an irrational number. Because by the
theorem 2,
$$lim_{n\rightarrow\infty}a_n\frac{c_{n+1}}{c_{n}}=lim_{n\rightarrow\infty}\frac{n!}{(n+1)!}=lim_{n\rightarrow\infty}\frac{1}{n+1}=0$$
Where $a_{n}=1$

{\bf Example 2.}
$$\theta=\sum_{n=1}^{\infty}\frac{n^4}{(n!)^5}=1+\frac{2^4}{(2!)^5}+\cdots+\frac{n^4}{(n!)^5}+\cdots
$$ is an irrational number, because by the theorem 2,$$a_n\frac{c_{n+1}}{c_{n}}=a_{n+1}\frac{b_{n}}{b_{n+1}}=\frac{1}{n+1}\rightarrow 0
(n\rightarrow \infty)$$ where $a_n=n^4, b_n=(n!)^5$.

Let's look at a more complicated example as follows:

{\bf Example 3.} Suppose that $r\geq 1$ is an integer, then
$sin\frac{1}{r}$ is an irrational number.

Since
$$sin\frac{1}{r}=\sum_{n=1}^{\infty}(-1)^{n-1}\frac{1}{(2n-1)!r^{2n-1}}=\sum_{n=1}^{\infty}c_n=\sum_{n=1}^{\infty}\frac{a_n}{b_n}$$
where $a_n=(-1)^{n-1},b_n=(2n-1)!r^{2n-1}$, then
$$lim_{n\rightarrow\infty}a_n\frac{c_{n+1}}{c_{n}}=lim_{n\rightarrow\infty}a_n\frac{b_{n}}{b_{n+1}}=lim_{n\rightarrow\infty}
\frac{(-1)^{n}(2n-1)!r^{2n-1}}{(2n+1)!r^{2n+1}}=0$$.
 The remaining
is to verify $\sum_{n=N}^{\infty}c_n\neq 0$
$$\sum_{n=N}^{\infty}c_n=\frac{(-1)^{N-1}}{(2N-1)!r^{2N-1}}+\frac{(-1)^{N}}{(2N+1)!r^{2N+1}}+\frac{(-1)^{N+1}}{(2N+3)!r^{2N+3}}+\cdots$$
When $N$ is odd
$$\sum_{n=N}^{\infty}c_n=(\frac{1}{(2N-1)!r^{2N-1}}-\frac{1}{(2N+1)!r^{2N+1}})+(\frac{1}{(2N+3)!r^{2N+3}}-\frac{1}{(2N+5)!r^{2N+5}})+\cdots$$
$$=\frac{(2N+1)!r^{2}-(2N-1)!}{(2N-1)!(2N+1)!r^{2N+1}}+\frac{(2N+5)!r^{2}-(2N+3)!}{(2N+3)!(2N+5)!r^{2N+5}}+\cdots >0$$
Similarly,when $N$ is even
$$\sum_{n=N}^{\infty}c_n=-(\frac{1}{(2N-1)!r^{2N-1}}-\frac{1}{(2N+1)!r^{2N+1}})-(\frac{1}{(2N+3)!r^{2N+3}}-\frac{1}{(2N+5)!r^{2N+5}})-\cdots<0$$
Thus for any natural number $N,\sum_{n=N}^{\infty}c_n\neq 0$,by the
theorem 1,$sin\frac{1}{r}$ is an irrational number.
 Since the sum of
two irrational numbers is not necessarily an irrational number, the
following theorem is interesting.
\begin{theorem}
Assume that $\alpha=\sum_{n=1}^{\infty}c_n$,where
$c_{n}=\frac{a_n}{b_n}>0$ and $\beta=\sum_{n=1}^{\infty}c_n^{'}$,
where$ c_n^{'}=\frac{a_n^{'}}{b_n^{'}}>0$ and above two number are
irrational numbers determined by the theorem 2 and satisfy the
following conditions:
$lim_{n\rightarrow\infty}\frac{a_{n+1}b_{n}b_{n}^{'}}{b_{n+1}}=0$
and
$lim_{n\rightarrow\infty}\frac{a_{n+1}^{'}b_{n}^{'}b_{n}}{b_{n+1}^{'}}=0$
then $\alpha +\beta$ is also an irrational number.
\end{theorem}
\begin{proof}
Let $\gamma=\alpha+\beta=\sum_{n=1}^{\infty}d_n$ and
$d_n=c_n+c_n^{'}=\frac{a_n}{b_n}+\frac{a_n^{'}}{b_n^{'}}=\frac{a_{n}b_n^{'}+b_{n}a_n^{'}}{b_{n}b_n^{'}}=\frac{\tilde{a_n
}}{\tilde{b_n}}$ where $\tilde{a_n }=a_{n}b_n^{'}+b_{n}a_n^{'},
\tilde{b_n}=b_{n}b_n^{'}$ then
$$lim_{n\rightarrow\infty}\tilde{a_n}\frac{d_{n+1}}{d_{n}}=lim_{n\rightarrow\infty}\tilde{a_{n+1}}\frac{\tilde{b_{n}}}{\tilde{b_{n+1}} }$$
$$=lim_{n\rightarrow\infty}(a_{n+1}b_{n+1}^{'}+b_{n+1}a_{n+1}^{'})\frac{b_{n}b_n^{'}}{b_{n+1}b_{n+1}^{'}}
=lim_{n\rightarrow\infty}\frac{a_{n+1}b_{n}b_{n}^{'}}{b_{n+1}}+lim_{n\rightarrow\infty}\frac{a_{n+1}^{'}b_{n}^{'}b_{n}}{b_{n+1}^{'}}=0$$
In additional,we notice that  $\tilde{b_n}|\tilde{b_{n+1}}$, so by
the theorem 2 we get $\alpha+\beta$ is an irrational number
\end{proof}
\newline
{\bf Example 4.}Let
$\alpha=\sum_{n=1}^{\infty}\frac{1}{2^{n!}},\beta=\sum_{n=1}^{\infty}\frac{1}{3^{n!}}$,we
can use above theorem to verify that $\alpha+\beta$ is an irrational
number as follows: First of all, $\alpha,\beta$ are irrational
numbers because of theorem 2,secondly, let
$b_n=2^{n!},b_n^{'}=3^{n!}$,then
$$\frac{b_{n}^{'}b_{n}}{b_{n+1}^{'}}=\frac{3^{n!}}{2^{nn!}}\rightarrow 0,(n\rightarrow \infty)$$
$$\frac{b_{n}^{'}b_{n}}{b_{n+1}^{'}}=\frac{2^{n!}}{3^{nn!}}\rightarrow 0,(n\rightarrow \infty)$$
so $\alpha+\beta$ is an irrational number.

Essentially, we can replace condition 1 of theorem 2 by a more
general condition as follows:
\begin{theorem}
Let series $\sum_{n=1}^{\infty}\frac{a_n}{b_n}$ where
$\frac{a_n}{b_n}>0$ are rational numbers and
$lim_{n\rightarrow\infty}\frac{a_{n+1}}{b_{n+1}}[b_1,b_2,\cdots
b_n]=0$ where $[b_1,b_2,\cdots b_n]$ denotes least common multiple
of $b_1,\cdots b_n$ then $\sum_{n=1}^{\infty}\frac{a_n}{b_n}$
converges an irrational number
\end{theorem}
\begin{proof}
Let $c_n=\frac{a_n}{b_n}$ and
$$c_n=\frac{a_{n}\frac{[b_1,b_2,\cdots b_n]}{b_n}}{[b_1,b_2,\cdots
b_n]}=\frac{\tilde{a_n}}{\tilde{b_n}}=\tilde{c_n}$$ where
$\tilde{a_n}=\frac{[b_1,b_2,\cdots
b_n]}{b_n},\tilde{b_n}=[b_1,b_2,\cdots b_n]$ and
$$lim_{n\rightarrow\infty}\tilde{a_n}\frac{\tilde{c_{n+1}}}{\tilde{c_{n}}}=
lim_{n\rightarrow\infty}\frac{\tilde{a_{n+1}}}{\tilde{b_{n+1}}}\tilde{b_n}=
lim_{n\rightarrow\infty}\frac{a_{n+1}}{b_{n+1}}[b_1,b_2,\cdots
b_n]=0$$ Obviously,$\tilde{b_n}|\tilde{b_{n+1}},n=1,2\cdots$ and the
series $\sum_{n=1}^{\infty}\tilde{c_{n}}$ satisfies conditions of
theorem 2 ,and
$$\sum_{n=1}^{\infty}\frac{a_n}{b_n}=\sum_{n=1}^{\infty}c_{n}=\sum_{n=1}^{\infty}\tilde{c_n}$$then we
finish the proof.
\end{proof}
\newline
{\bf Example 5.}
$$\theta=\frac{1}{p_{2^{2^{1!}}}}+\frac{1}{p_{2^{2^{2!}}}}+\cdots+\frac{1}{p_{2^{2^{n!}}}}+\cdots$$
is a irrational number, where $p_{2^{2^{n!}}}$ is $2^{2^{n!}}$-th
prime number. Let's show it as follows:
\newline
We need the famous result\cite{Hua}: let $p_n$ is $n$-th prime
number, there exists two positive numbers such that
$c_{1}n\ln{n}<p_n<c_{2}n\ln{n}$,then
$$\frac{a_{n+1}}{b_{n+1}}[b_1,b_2,\cdots
b_n]=\frac{p_{2^{2^{1!}}}p_{2^{2^{2!}}}\cdots
p_{2^{2^{n!}}}}{p_{2^{2^{(n+1)!}}}}$$
$$<\frac{c_{2}2^{2^{1!}}\ln{2^{2^{1!}}}c_{2}2^{2^{2!}}\ln{2^{2^{2!}}}\cdots c_{2}2^{2^{n!}}\ln{2^{2^{n!}}}}
{c_{1}2^{2^{(n+1)!}}\ln{2^{2^{(n+1)!}}}}$$
$$=\frac{2^{2^{1!}+2^{2!}+\cdots+2^{n!}+n\log_{2}C}}{2^{2^{(n+1)!}}}
\frac{2^{1!+2!+\cdots+n!}}{2^{(n+1)!}}k(\ln2)^{n-1}$$ Where
$k=\frac{1}{c_1},c=c_{2}$ Since
$$1!+2!+\cdots+n!\leq nn!<(n+1)!$$, $\frac{2^{1!+2!+\cdots+n!}}{2^{(n+1)!}}<1 .$
Also, there exists $N$ such that when $n\geq N$ $\log_{}2C<2^{n!}$.
Therefore
$$2^{2^{1!}+2^{2!}+\cdots+2^{n!}+n\log_{2}C}\leq n2^{n!}+n2^{n!} \leq 2n2^{n!}\leq 2^{(n+1)!}$$
we get
$$\frac{2^{2^{1!}+2^{2!}+\cdots+2^{n!}+n\log_{2}C}}{2^{2^{(n+1)!}}}\leq 1$$
Thus when $n\geq N$
$$lim_{n\rightarrow\infty}\frac{a_{n+1}}{b_{n+1}}[b_1,b_2,\cdots
b_n]\leq k(\ln2)^{n-1}\rightarrow 0(n\rightarrow \infty)$$ ,by
theorem 4,$\theta $ is an irrational number
\newline
The following theorem shows that  condition 2 of theorem 2 is also
necessary in some special case.

\begin{theorem}
Assume that the sequence $c_{n} =\frac{1}{a^{P_{m}(n)}}$ ,where
$a\geq 2$ is an integer and
$P_{m}(x)=b_{0}x^{m}+b_{1}x^{m-1}+\cdots+b_{m-1}x+b_m$ is an
polynomial with positive integer coefficients,then series
$\sum_{n=1}^{\infty}c_n$ converges to an irrational number if and
only if $lim_{n\rightarrow\infty}\frac{c_{n+1}}{c_{n}}=0$
\end{theorem}
\begin{proof}
Obviously, $P_{m}(n)<P_{m}(n+1)$,so
$a^{P_{m}(n)}|a^{P_{m}(n+1)},n=1,2,\cdots$ which satisfies condition
1 of  theorem 2 and

\begin{eqnarray}
P_{m}(n+1)-P_{m}(n) &=& b_{0}(n+1)^{m}+\cdots
+b_{m}-(b_{0}^{m}+\cdots+b_{m}
+b_m)\nonumber\\
                    &=&b_{0}(n^{m}+mn^{m-1}+\cdots)+\cdots +b_m-(b_{0}^{m}+\cdots+b_m
+b_m)\nonumber\\
                    &= &mb_{0}n^{m-1}+l_{1}n^{m-2}+l_{2}n^{m-3}+\cdots
\end{eqnarray}
then
$$\frac{c_{n+1}}{c_n}=\frac{a^{P_{m}(n)}}{a^{P_{m}(n+1)}}
=\frac{1}{a^{P_{m}(n+1)-P_{m}(n)}}=\frac{1}{a^{mb_{0}n^{m-1}+l_{1}n^{m-2}+
\cdots}}$$ Thus
\begin{eqnarray}
lim_{n\rightarrow\infty}\frac{c_{n+1}}{c_{n}}=\left \{
\begin{array}{ll}
0 & m\geq 2 \\
a^{b_0} & m=1
\end{array}
\right.
\end{eqnarray}
By theorem 2, when $m\geq 2$ series $\sum_{m=1}^{\infty}c_n$
converges to an irrational number . and
$lim_{n\rightarrow\infty}\frac{c_{n+1}}{c_{n}}\neq 0$ means $m=1$
then
$$\sum_{n=1}^{\infty}c_n=\sum_{n=1}^{\infty}\frac{1}{a^{b_{0}n+b_{1}}}=\frac{1}{a^{b_1}(a^{b_0}-1)}$$
is a rational number.
\end{proof}
\begin{theorem}
Let $\theta=\sum_{n=1}^{\infty}\frac{a_n}{b_n}$,where
$\frac{a_n}{b_n}>0 (n=1,2\cdots)$ are rational number, and assume
that $ f(b_n)$ is a function of $b_n$  and $ f(b_n)>0$ ,furthermore
 if  the following conditions are satisfied:
\newline
(1)$b_1<b_2<\cdots$ and $b_n|b_{n+1},n=1,2,\cdots$ \newline
(2)$f(b_n)>0$ and $\frac{f(b_{n+1})}{f(b_n)}<\frac{1}{2}$ ($n$ is
big enough)\newline (3)$\frac{f(b_{n})}{b_{n+1}}a_{n+1}<\frac{1}{2}$
($n$ is big enough)\newline (4)$\frac{b_n}{f(b_n)}\rightarrow 0
(n\rightarrow \infty)$
\newline
 then
\newline
(1)$\theta$ is an irrational number
\newline (2) When $n$ is big enough, there exists infinite fractions
$\frac{c_n}{b_n}$ such that $|\theta -
\frac{c_n}{b_n}|<\frac{1}{f(b_n)}$
\end{theorem}
\begin{proof}
According condition 3,when $n$ is big enough, we have $\frac
{f(b_{n})}{b_{n+1}a}a_{n+1}<\frac{1}{2}$ or equivalently
$\frac{a_{n+1}b_n}{b_{n+1}}<\frac{b_n}{2f(b_n)}$,so
$$lim_{n\rightarrow\infty}a_n\frac{c_{n+1}}{c_{n}}=lim_{n\rightarrow\infty}\frac{a_{n+1}b_n}{b_{n+1}}=0$$
In the last step we use condition 4.By theorem 2, $\theta$ is an
irrational number.
\newline Let's prove the second part.

$$|\theta-\frac{a_1}{b_1}-\frac{a_2}{b_2}-\cdots-\frac{a_n}{b_n}|=\frac{a_{n+1}}{b_{n+1}}+\frac{a_{n+2}}{b_{n+2}}+\frac{a_{n+3}}{b_{n+3}}+\cdots$$

$$=\frac{1}{f(b_n)}(\frac{f(b_{n})a_{n+1}}{b_{n+1}}+\frac{f(b_{n})a_{n+2}}{b_{n+2}}+\frac{f(b_{n})a_{n+3}}{b_{n+3}}+\cdots)$$

$$=\frac{1}{f(b_n)}(\frac{f(b_{n})a_{n+1}}{b_{n+1}}+\frac{f(b_{n})}{f(b_{n+1})}\frac{f(b_{n+1})a_{n+2}}{f(b_{n+2})}+\frac{f(b_{n})}{f(b_{n+1})}\frac{f(b_{n+1})}{f(b_{n+2})}\frac{f(b_{n+2})a_{n+3}}{f(b_{n+3})}+\cdots)$$

$$<\frac{1}{f(b_n)}(\frac{1}{2}+\frac{1}{2^2}+\frac{1}{2^3}+\cdots)=\frac{1}{f(b_n)}$$

We use condition 2 and 3 in the last two steps.

Let
$c_{n}=b_{n}(\frac{a_1}{b_1}+\frac{a_2}{b_2}+\cdots+\frac{a_n}{b_n})$,since
$\frac{a_n}{b_n}>0,(n=1,2,\cdots)$,
$\frac{c_n}{b_n}<\frac{c_{n+1}}{b_{n+1}}$ Thus there are infinite
number of $\frac{c_n}{b_n}$ satisfy
$$|\theta-\frac{c_n}{b_n}|=|\theta-\frac{a_1}{b_1}-\frac{a_2}{b_2}-\cdots-\frac{a_n}{b_n}|<\frac{1}{f(b_n)}$$
(when $n$ is big enough )That proves the theorem
\end{proof}
\newline
{\bf Remark.} In above theorem, condition 2 and 3 can be replaced by
$lim_{n\rightarrow\infty}\frac{f(b_n)}{f(b_{n+1})}=l<\frac{1}{2}$
and
$lim_{n\rightarrow\infty}\frac{f(b_{n})}{b_{n+1}}a_{n+1}=k<\frac{1}{2}$
Using theorem 6 and two following known results, we can get two
useful theorems , one is about how to determine a given number is
transcendental number, the other is about complex analytic dynamics.

{\bf Theorem(K.Roth).} Let $\theta$ be a $n\geq 2$ degree algebraic
number, then for any given $\epsilon > 0$, there exists only finite
positive integer pairs $x,y$ such that $|\theta-\frac{x}{y}|<
\frac{1}{y^{2+\epsilon}}$

{\bf Theorem (H.Cremer)\cite{Cremer}} If irrational number $\theta$
satisfies the condition that there exists infinite positive integers
such that $|\theta-\frac{n}{m}|\leq \frac{1}{m^{d^{m}-1}}$,
indifferent fixed point $z=0$ of polynomial
$f(z)=z^{d}+\cdots+e^{2\pi i \alpha}$ belongs to Julia set. First of
all, we use theorem 6 and Roth theorem to derive following theorem:
\begin{theorem}
Let $\theta=\sum_{n=1}^{\infty}\frac{a_n}{b_n}$,where
$\frac{a_n}{b_n}>0 (n=1,2\cdots)$ are rational numbers which satisfy
\newline
(1)$b_1<b_2<\cdots$ and $b_n|b_{n+1},n=1,2,\cdots$
\newline
(2) for some $\epsilon >0,$
$\frac{a_{n+1}b_{n}^{2+\epsilon}}{b_{n+1}}<\frac{1}{2}$ ($n$ is big
enough)
\newline
then $\theta=\sum_{n=1}^{\infty}\frac{a_n}{b_n}$ is a transcendental
number.
\end{theorem}
\begin{proof}
 In the theorem 6,let's take $f(b_n)=b_{n}^{2+\epsilon}$,then it's easy to
verify $f(b_n),n=1,2,\cdots$ satisfy  condition 1 and 3 of theorem
6, we only need to check condition 2 and 4.
\newline
Since $\epsilon>0$ and $a_n,b_n$ are positive integers, we have
$$\frac{f(b_{n})}{f(b_{n+1})}=\frac{b_{n}^{2+\epsilon}}{b_{n+1}^{2+\epsilon}}< \frac{b_{n}^{2+\epsilon}}{b_{n+1}}<\frac{a_{n+1}b_{n}^{2+\epsilon}}{b_{n+1}}<\frac{1}{2}$$ ,thus we pass condition 2 of theorem 6.
Also because $b_1<b_2<\cdots$ and
$\epsilon>0$,$$\frac{b_n}{f(b_n)}=\frac{b_n}{b_{n}^{2+\epsilon}}=\frac{1}{b_{n}^{1+\epsilon}}\rightarrow
0(n\rightarrow \infty)$$, then we finish checking all conditions of
theorem 6 get satisfied. By theorem 6,
 $\theta$ is an irrational
number, or equivalently, it's not a first order algebraic number,by
the conclusion 2, when $n$ is big enough, there exists infinite
fractions $\frac{c_n}{b_n}$ satisfy $|\theta-\frac{c_n}{b_n}|<
\frac{1}{{b_n}^{2+\epsilon}}$,by Roth theorem we get $\theta$ is a
transcendental number.
\end{proof}
\newline
{\bf Example 6.} $\sum_{m=1}^{\infty}\frac{1}{10^{m!}}$ is an
transcendental number.
\newline
Because by taking $a_{n}=1,n=1,2,\cdots, b_{n}=10^{n!}$ and
$\epsilon=1$,
$\frac{a_{n+1}b_{n}^{2+\epsilon}}{b_{n+1}}=\frac{(10^{n!})^3}{10^{{n+1}!}}=\frac{1}{10^{n!(n-2)}}\rightarrow
0 (n\rightarrow \infty)$,by theorem 7,
$\sum_{m=1}^{\infty}\frac{1}{10^{m!}}$ is an transcendental number.
\newline
{\bf Example 7.} $\sum_{n=1}^{\infty}\frac{3^n}{2^{3^n}}$ is an
transcendental number.
\newline
Because by taking $a_{n}=3^n,n=1,2,\cdots, b_{n}=2^{3^n}$ and
$\epsilon=\frac{2}{3}$,
$\frac{a_{n+1}b_{n}^{2+\epsilon}}{b_{n+1}}=9\frac{3^{n-1}}{2^{3^{n-1}}}
\rightarrow 0 (n\rightarrow \infty)$,by theorem 7,
$\sum_{n=1}^{\infty}\frac{3^n}{2^{3^n}}$is an transcendental
number.
\begin{theorem}
Let $\theta=\sum_{n=1}^{\infty}\frac{a_n}{b_n}$,where
$\frac{a_n}{b_n}>0 (n=1,2\cdots)$ are rational numbers which satisfy
\newline
(1)$b_1<b_2<\cdots$ and $b_n|b_{n+1},n=1,2,\cdots$
\newline (2)
$\frac{a_{n+1}b_{n}^{d^{b_n}-1}}{b_{n+1}}<\frac{1}{2}$ ($d\geq 2$ is
an integer and $n$ is big enough)
\newline
Then indifferent fixed point $z=0$ of polynomial
$f(z)=z^{d}+\cdots+e^{2\pi i \alpha}$ belongs to Julia set
\end{theorem}
\begin{proof}
Set $f(b_n)=b_{n}^{d^{b_n}-1}$,it's easy to verify $f(b_n)$ satisfy
condition 1 and 3 of theorem 6. We only need to verify condition 2
and 4. Since $b_n|b_{n+1},(n=1,2,\cdots)$ and $b_{n+1}\geq 2b_{n}$
,we have
$$\frac{f(b_{n})}{f(b_{n+1})}=\frac{b_{n}^{d^{b_n}-1}}{b_{n+1}^{d^{b_{n+1}}-1}}\leq
\frac{(\frac{1}{2}b_{n+1})^{d^{\frac{1}{2}b_{n+1}}-1}}{b_{n+1}^{d^{b_{n+1}}-1}}=
\frac{(\frac{1}{2})^{d^{\frac{1}{2}b_{n+1}}-1}b_{n+1}^{d^{\frac{1}{2}b_{n+1}}-1}}{b_{n+1}^{d^{b_{n+1}}-1}}$$

$$<(\frac{1}{2})^{d^{\frac{1}{2}b_{n+1}}-1}<\frac{1}{2}$$
thus condition (2) is satisfied. Let's verify condition (4), Since
when $n\geq 2$,$b_n \geq 2$ and $d\geq 2$,$d^{b_n}-1\geq 3$, thus
when $n\geq 2$
$$\frac{b_n}{f(b_n)}=\frac{b_n}{b_{n}^{d^{b_n}-1}}\leq \frac{1}{b_n^2 }\rightarrow 0(n\rightarrow \infty)$$
By theorem 6,$\theta$ is an irrational number and there exists
infinite fractions $\frac{c_n}{b_n}$ such that
$|\theta-\frac{c_n}{b_n}|<\frac{1}{b_{n}^{d^{b_n}-1}}$ ($n$ is big
enough), then by Cremer theorem, we get our result
\end{proof}
\newline
{\bf Example 8.} In order to illustrate this example, we need some
notation to describe a special series so called "nth exponential
floor"as follows:
\newline

            Set $[a_n,a_{n-1},\cdots, a_1]_{n} =f_n$ where $f_n$
            is defined inductively by $f_{k+1}=(a_{k+1})^{f_k},k=1,2\cdots, n-1$
            and $f_{1}=a_1$

           For any positive integer $d\geq2$, let $b_{n}=[d,\cdots,d,nd]_{2n}$
           and $\theta=\sum_{n=1}^{\infty}\frac{1}{b_n}$, we will
           show that indifferent fixed point $z=0$ of polynomial $g(z)=Z^{d}+\cdots+e^{2\pi
           i\theta}z$  belongs to Julia set.
           Let's check $b_n$ satisfy conditions of theorem 8,
           condition 1 is obvious and by noticing
           $[d,\cdots,d,nd]_{2n}=d^{[d,\cdots,d,nd]_{2n-1}}$,we have
           $$\frac{a_{n+1}b_{n}^{d^{b_n}-1}}{b_{n+1}}=\frac{([d,\cdots,d,nd]_{2n})^{[d,\cdots,d,nd]_{2n+1}-1}}{[d,\cdots,d,(n+1)d]_{2(n+1)}}$$
           $$<\frac{([d,\cdots,d,nd]_{2n})^{[d,\cdots,d,nd]_{2(n+1)}}}{[d,\cdots,d,(n+1)d]_{2(n+1)}}$$
           $$=\frac{d^{([d,\cdots,d,nd]_{2n-1})([d,\cdots,d,nd]_{2n+1})}}{[d,\cdots,d,(n+1)d]_{2(n+1)}}
           =\frac{d^{d^{[d,\cdots,d,nd]_{2n-2}+[d,\cdots,d,nd]_{2n}}}}{[d,\cdots,d,(n+1)d]_{2(n+1)}}$$
           $$=\frac{d^{d^{[d,\cdots,d,nd]_{2n-2}+[d,\cdots,d,nd]_{2n}}}}{d^{d^{[d,\cdots,d,(n+1)d]_{2n}}}}$$

        We notice that  $\frac{[d,\cdots,d,nd]_{2n-2}+[d,\cdots,d,nd]_{2n}}{[d,\cdots,d,(n+1)d]_{2(n)}}\rightarrow 0(n\rightarrow\infty)$

         That means when $n$ is big enough,$\frac{a_{n+1}b_{n}^{d^{b_n}-1}}{b_{n+1}}<\frac{1}{2}$ and we finish checking this example satisfies all
         conditions of theorem 8 and thus indifferent fixed point of $g(z)$
         belongs to Julia set.
\newline
At the end of paper, we are going to propose a conjecture which
relates to theorem 2. To do this, we need following definition
firstly.
\begin{definition}
Let $\alpha$ be an irrational number ,if $\alpha$ satisfies
following conditions:
\newline
(1) $\alpha=\sum_{n=1}^{\infty}c_{n}, $ where
$c_{n}=\frac{a_{n}}{b_{n}} ,(n=1,2,\cdots)$ and $a_n,b_n$ are
positive integers.
\newline
(2)  $$b_n |b_{n+1},(n=1,2,\cdots)$$
\newline
(3)  $$lim_{n\rightarrow\infty}a_n\frac{c_{n+1}}{c_{n}}=0$$

We call the irrational number $\alpha$ has $E$ rational
approximation.
\end{definition}

{\bf Conjecture.}Every positive irrational number has $E$ rational
approximation.

{\bf Remark.}  The positive answer of above conjecture will give an
explicit character  of any positive irrational number.

\end{document}